\newif\ifsmfart
\IfFileExists{smfart.cls}
   {\documentclass[12pt,english]{smfart}\smfarttrue}
   {\message{^^J*** This file should be typeset with smfart.cls ***^^J^^J}
    \documentclass[12pt]{amsart}}

\usepackage{amssymb}

\IfFileExists{mathrsfs.sty}
  {\usepackage{mathrsfs}\let\mathcal\mathscr}
  {\let\mathscr\mathcal}

\advance\headheight 2pt
\calclayout
\allowdisplaybreaks[2]

\numberwithin{equation}{section}
\makeatletter

\let\c@equation\c@subsection
\let\cl@equation\cl@subsection
\def\theequation{\thesection.\arabic{equation}}

\def\l@table{\@tocline{0}{3pt plus2pt}{0pt}{}{\itshape}}

\makeatother

\theoremstyle{plain}
\newtheorem{prop}[subsection]{Proposition}
\newtheorem{thm}[subsection]{Theorem}
\newtheorem{cor}[subsection]{Corollary}

\newtheorem{lem}[subsection]{Lemma}

\theoremstyle{definition}

\theoremstyle{remark}

\def\A{{\mathbf A}}
\def\C{{\mathbf C}}
\def\F{{\mathbf F}}
\def\G{{\mathbf G}}

\def\P{{\mathbf P}}
\def\Q{{\mathbf Q}}
\def\R{{\mathbf R}}

\def\Z{{\mathbf Z}}

\let\ra\rightarrow

\let\epsilon\varepsilon \let\eps\epsilon
\let\epsilon\varepsilon
\let\phi\varphi

\let\leq\leqslant
\let\geq\geqslant

\def\eff{{\text{\upshape eff}}}
\def\abs#1{\left\lvert{#1}\right\rvert}
\def\norm#1{\left\|{#1}\right\|}

\def\DeclareMathOperator#1#2{\def #1{\operatorname{#2}}}
\DeclareMathOperator{\Re}{Re} 
\DeclareMathOperator{\Pic}{Pic}

\DeclareMathOperator{\Spec}{Spec}
\DeclareMathOperator{\Proj}{Proj}
\DeclareMathOperator{\rk}{rk}

\DeclareMathOperator{\Val}{Val}
\DeclareMathOperator{\vol}{vol}

\def\dx{\,{\mathrm d\mathbf x}}

\title[Points of bounded height, II]
      {Points of bounded height \\
       on equivariant compactifications \\ of vector groups, II}

\author{Antoine Chambert-Loir}
\address{
Institut de math\'ematiques de Jussieu, Boite 247 \\ 
4 place Jussieu \\ 
F-75252 Paris Cedex 05 }
\email{chambert@math.jussieu.fr}

\author{Yuri Tschinkel}
\address{Dept. of Mathematics, U.I.C.\\
Chicago, (IL) 60607-7045,  U.S.A. }
\email{yuri@math.uic.edu}

\begin{document}
\def\smfandname{\&}
\date{Submitted on the xxx archive on  April 2nd, 1999}

\begin{abstract}
We prove asymptotic formulas  for the number
of rational points of bounded height
on certain equivariant compactifications of the affine space.
\end{abstract}
 
\ifsmfart 
 \begin{altabstract}
 Nous \'etablissons un d\'eveloppement asymptotique
 du nombre de points rationnels de hauteur born\'ee sur
 certaines compactifications \'equi\-vari\-antes de l'espace affine.
 \end{altabstract}
\fi

\maketitle

\tableofcontents

\section*{Introduction}
\def\thesubsection{\arabic{subsection}}
In the last decade, there has been much interest in establishing
asymptotics for the number of points of bounded height on algebraic
varieties defined over number fields. 
Yu.~Manin and V.~Batyrev~\cite{batyrev-m90} have
formulated conjectures describing such asymptotics in geometrical
terms. These conjectures have been further refined by 
E.~Peyre in \cite{peyre95}.

More precisely, let $X$ be a smooth projective algebraic variety defined
over a number field $F$ and $H:X(F)\ra\R_{>0}$ an exponential
height function on the set of rational points of $X$ defined by some 
metrized ample line bundle~$\mathscr L$.
One wants to relate the asymptotic behaviour of the counting function
\[ N(U,\mathscr L,B) = \# \{ x\in U(F)\,;\, H(x)\leq B\} \]
to geometric invariants of $X$, such as the cone of effective line
bundles and the (anti)-canonical line bundle of $X$. Here, $U$ is
a sufficiently small Zariski dense open subset; its presence is made
necessary by possible ``accumulating subvarieties'', which
contain more rational points than their complement in $X$. 
If $X$ is a Fano variety and $\mathscr L=K_X^{-1}$, one expects that
\[ N(U,K_X^{-1},B) \sim \frac{\Theta(X)}{(r-1)!} B (\log B)^{r-1} \]
where $r=\rk \Pic(X)$ and $\Theta(X)$ is the product
of three numbers:  a Tamagawa constant which measures
the volume of the closure of rational points in the adelic points 
$\overline{X(F)}\subset X(\A_F)$ 
with respect to the metrization,
a rational number defined in terms of the cone of effective divisors
and the order of the non-trivial part of the Brauer group of $X$.

Such a description cannot hold universally (see the example by
V.~Batyrev and Yu.~Tschinkel~\cite{batyrev-t96}), but there are
two classes of algebraic varieties where it does hold:
those for which the \emph{circle method} in analytic
number theory applies, and those possessing many symmetries, 
such as an \emph{action} (with a dense orbit) of a linear algebraic group.
The circle method is concerned with 
complete intersections of small degree and
small codimension in projective space. 
They have moduli, but only few projective embeddings;
the Picard group is $\Z$. 
As a reference, let us mention the papers by B.~Birch~\cite{birch62}
and W.~Schmidt~\cite{schmidt85}. 
The other approach leads, via harmonic analysis
on the adelic points of the corresponding group, to a proof of 
conjectured asymptotic formulas for toric varieties (see~\cite{batyrev-t98b})
or for generalized flag varieties (using Langlands' work on
Eisenstein's series, see~\cite{franke-m-t89}).
These have Picard groups of higher ranks, but 
no deformations due to the \emph{rigidity} of reductive groups.

In this paper we treat certain equivariant compactifications 
of vector groups. In a previous
paper~\cite{chambert-loir-t99c},
we had established asymptotic formulas for 
blow-ups of $\P^2$ in any number of points \emph{on a line}. 
Here we work out the case of blow-ups of a projective
space $\P^n$ of dimension at least $3$ in a smooth codimension 2
subvariety contained in a hyperplane. It should be clear to the reader that
these varieties admit deformations (they are
parametrized by an open subset of an appropriate Hilbert scheme).

More precisely, 
let $f\in\Z[x_1,\dots,x_n]$ be a homogeneous polynomial of degree~$d$
and $X\rightarrow \P^n=\Proj(\Z[x_0,\dots,x_n])$ be the
blow-up of the ideal generated by $(x_0,f)$.
Suppose that the hypersurface defined by $f$ in $\P^{n-1}_\C$
is smooth and let $U\simeq\A^n$ be the inverse image in $X$
of $\A^n\subset\P^n$.
Then, $X_\C$ is a smooth projective variety, with Picard group $\Z^2$
and trivial Brauer group. Moreover, $X_\C$ is an equivariant
compactification of $\G_a^n$.
There is a natural metrization on $K_X^{-1}$ (recalled
below) which allows to define the height function and 
the \emph{height zeta function}
\[ Z(U,K_X^{-1},s) = \sum_{x\in U(\Q)} H_{K_X^{-1}}(x)^{-s},  \]
The series converges absolutely for $\Re(s)\gg 0$.
Our main theorem is:
\begin{thm}
There exists a function $h$ which is holomorphic in the domain
$\Re(s)> 1-\frac 1n$ such that
\[ Z(U,K_X^{-1},s) = \frac{h(s)}{(s-1)^2}\quad\text{and}\quad
      h(1)=\Theta(X)\neq 0. \]
\end{thm}

A standard Tauberian theorem implies that $X$ satisfies Peyre's
refinement of Manin's conjecture:
\begin{cor}
We have the following asymptotic formula:
\[ N(U,K_X^{-1},B) \sim {\Theta(X)} B \log(B)\]
as $B$ tends to infinity.
\end{cor}

In fact, we will prove asymptotics for every $\mathscr L$ on $X$
such that its class is contained in the interior of the effective
cone $\Lambda_\eff(X)$. Moreover, we will prove estimates for the
growth of $Z(s)$ in vertical strips in the neighbourhood of
$\Re(s)=1$. It is well known that this implies a more precise asymptotic
expansion for the counting function $N(U,\mathscr L,B)$,
see Theorem~\ref{thm:peyre} and its corollary at the end of the paper.

\newpage

\section{Geometry, heights}
\def\thesubsection{\thesection.\arabic{subsection}}

Let $f\in\Z[x_1,\dots,x_n]$ be a homogeneous polynomial of degree~$d$
with coprime coefficients and $\pi:X\ra\P^n$ the blow-up of
the ideal $(x_0,f)$ in $\P^n=\Proj(\Z[x_0,\dots,x_n])$.
We denote by $Z_f$ the hypersurface defined by $f$ in $\P^{n-1}$.
Throughout the paper, we assume that $Z_{f,\C}$ is smooth,
irreducible and  that it doesn't contain any hyperplane. In other words,
$n\geq 3$ and $d\geq 2$.
The universal property of blowing up implies
that the scheme $X$ is an equivariant compactification
of the additive group $\G_a^n=\Spec(\Z[x_1,\dots,x_n])$.

Denote by $D_1$ the exceptional divisor in $X$ and by $D_0$
the strict transform of the divisor $x_0=0$ in $\P^n$.
Let $U\simeq\G_a^n$ be the inverse image of $\G_a^n$ under $\pi$.
We identify rational points in $U$ with their image in the
affine space $\G_a^n\subset\P^n$.

If $\mathbf s\in\C^2$, denote $D(\mathbf s)=s_0
[D_0]+s_1[D_1]\in\Pic(X)\otimes_\Z{\C}$.

The following proposition summarizes the geometric facts needed in the
sequel.

\begin{prop}
The classes of the divisors $D_0$ and $D_1$ form a basis of $ \Pic (X)$.
For $\mathbf s=(s_0,s_1)\in\Z^2$,
the divisor class $D(\mathbf s)$ is effective iff $s_0\geq 0$ and $s_1\geq 0$.
The variety $X_\Q$ is smooth; its anticanonical line bundle has class
$D(n+1,n)$.
\end{prop}
\begin{proof}
See~\cite{chambert-loir-t99c}, Prop.~1.3 and Prop.~1.6
or~\cite{hartshorne77}, chap.~II, \S\,8.
\end{proof}

We now define height functions on $X$. We denote
by $\Val(\Q)=\{2,3,\dots,\infty\}$ the set of places of $\Q$.
If $p$ is a prime number
and $\mathbf x\in\G_a^n(\Q_p)$, let
$\norm{\mathbf  x}_p=\max(\abs{x_1}_p,\dots,\abs{x_n}_p)$ and 
define the functions $H_{D_1,p}$ and $H_{D_0,p}$ by 
\begin{align}
H_{D_1,p}(\mathbf x)^{-1}
& = \max\big( \frac{1}{\max(1,\norm{\mathbf x}_p)},
    \frac{\abs{f(\mathbf x)}_p}{\max(1,\norm{\mathbf x}_p)^d} \big)  \\
H_{D_0,p}(\mathbf x)^{-1}
& = \frac{H_{D_1,p}(\mathbf x)}{\max(1,\norm{\mathbf x})}.
\end{align}
At the archimedian place of $\Q$, define the local
height functions
by replacing maximums by the square root of the sum of squares.
For any place $v$ of $\Q$ and any $\mathbf s=(s_0,s_1)\in\C^2$, we set 
\begin{equation}
H_v(\mathbf s;\mathbf x)=H_{D_0,v}(\mathbf x)^{s_0}
       H_{D_1,v}(\mathbf x)^{s_1}. 
\end{equation}
Finally, we define a global height pairing 
\begin{equation}
   H:\Pic(X)_\C \times \G_a^n(\A_\Q) \ra\C^*,
\quad
H(\mathbf s;\mathbf x) = \prod_{v\in\Val(\Q)} H_v(\mathbf s;\mathbf x_v). 
\end{equation}

\begin{prop}
If $\mathscr L\in\Pic(X)$,
the function $\mathbf x\mapsto H(\mathscr L;\mathbf x)$
on $\G_a^n(\Q)$ is an exponential height in the sense of Weil.
\end{prop}
\begin{proof}
See~\cite{chambert-loir-t99c}, (1.12), (1.13) and (2.2).
\end{proof}

The height zeta function is then defined by the series
\begin{equation}
Z(\mathbf s) = \sum_{\mathbf x\in \G_a^n(\Q)} H(\mathbf s;\mathbf x)^{-1}.
\end{equation}
It converges \emph{a priori}
for all $\mathbf s\in\C^2$ such that $D(\mathbf s)$ is sufficiently
ample, i.e.~if $\Re(s_0-s_1)$ and $\Re(s_0)$ are big enough.

Let $\psi=\prod_v\psi_v:\G_a(\A_\Q)\ra\C^*$ be the standard
additive character of $\A_\Q$. If $\mathbf a\in\Q^n$, we
define 
\[ \psi_{\mathbf a}(\mathbf x)=\psi(\langle\mathbf a,\mathbf x\rangle).\]
We use the standard self-dual Haar measure $\dx$
on $\G_a^n(\A_\Q)$.
For any $\mathbf a\in\Q^n$, define the Fourier transform
\[ \hat H(\mathbf s;\psi_{\mathbf a}) = \int_{\G_a^n(\A_\Q)} H(\mathbf
s;\mathbf x)^{-1}\, \dx. 
\]
It is the product of the local Fourier transforms $\hat H_v(\mathbf
s;\psi_{\mathbf a})$.

For $\mathbf s\in\C^2$ such that both sides converge absolutely,
we have the following identity:
\begin{equation}
Z(\mathbf s) = \sum_{\mathbf a\in\Z^n} \hat H(\mathbf s;\psi_{\mathbf a}).
\end{equation}
This is a consequence of the usual Poisson formula,
see~\cite{chambert-loir-t99c}, end of~\S\,2.

In the following sections 
we determine the domain of absolute convergence of the right hand side
and prove that $Z(\mathbf s)$  admits a meromorphic continuation beyond this
domain.

\section{The local Fourier transform at the trivial character}

We denote by $S$ the minimal set of primes
such that $Z_f\subset\P^{n-1}_\Z$
is smooth over $\Spec \Z[S^{-1}]$.
Let $p$ be a prime number.

\subsection{Decomposition of the domain}
We define subsets of $\Q_p^n$  as follows:
\begin{itemize}
\item $U(0)=\Z_p^n$;
\item if $0\leq \beta < \alpha$, $U_1(\alpha,\beta)$ is
the set of $\mathbf x\in\Q_p^n$ such that
$\norm{\mathbf x}=p^\alpha$ and $\abs{f(\mathbf x)}=p^{d\alpha-\beta}$;
\item if $\alpha\geq 1$, $U_1(\alpha)$ is the set of $\mathbf x\in\Q_p^n$
such that $\norm{\mathbf x}=p^\alpha$ and $\abs{f(\mathbf x)}\leq
p^{(d-1)\alpha}$;
\item if $\alpha\geq 1$, $U(\alpha)$ is the  set of $\mathbf x\in\Q_p^n$
such that $\norm{\mathbf x}=p^\alpha$ and $\abs{f(\mathbf x)}=p^{d\alpha}$.
\end{itemize}

The local height function is constant on each of these subsets.
Namely, if $\mathbf x\in U(0)$, $H_{D_0,p}=H_{D_1,p}=1$.
If $\mathbf x\in U_1(\alpha,\beta)$,
$H_{D_0,p}= p^{\alpha-\beta}$ and $H_{D_1,p}=p^\beta$.
On $U(\alpha)$, $H_{D_0,p}=p^\alpha$ and $H_{D_1,p}=1$.
Finally, if $\mathbf x\in U(\alpha)$, then 
$H_{D_0,p}=1$ and $H_{D_1,p}=p^\alpha$.

\subsection{Volumes}
Denote by 
\[ \tau_p(f) = \big( 1-\frac1p \big) \frac{\# Z_f(\F_p)}{p^{n-2}}.
\]
The Weil conjectures proved by Deligne imply that $\tau_p(f)=1+O(1/p)$.
In a much more elementary way,
it follows from Lemma~\ref{lem:weil} below that $\tau_p(f)$ is bounded
as $p$ varies.

\begin{lem}
For $p\not\in S$, we have
\addtocounter{equation}{-1}
\begin{subequations}
\begin{align}
\vol(U(0)) &= 1 \\
\vol (U_1(\alpha,\beta)) &= \frac{p-1}p \tau_p(f) p^{n\alpha -\beta} \\
\vol (U_1(\alpha)) &= \tau_p(f) p^{(n-1)\alpha} \\
\vol (U(\alpha)) & = \big(1-p^{-n} - p^{-1}\tau_p(f)\big)p^{n\alpha}.
\end{align}
\end{subequations}
\end{lem}

\begin{proof}
For $\beta\geq 1$, let $\Omega(\beta)$ be the set of $\mathbf x\in\Z_p^n$
such that $\norm{\mathbf x}=1$ and $\abs{f(\mathbf x)}\leq p^{-\beta}$.
By definition, 
\[ \vol(\Omega(\beta)) = p^{-n\beta} p^{\beta-1}(p-1) \#
Z_f(\Z/p^\beta\Z).\]
As $Z_f$ is smooth of pure dimension $n-2$ over $\Z_p$, Hensel's lemma implies
that 
\[ \# Z_f(\Z/p^\beta\Z)=p^{(\beta-1)(n-2)}\# Z_f(\F_p). \]
Consequently,
\[ \vol(\Omega(\beta)) = (p-1) p^{-\beta-1} \frac{\# Z_f(\F_p)}{p^{n-2}}
    = \tau_p(f)p^{-\beta} . \]
As $ U_1(\alpha)=p^{-\alpha} \Omega(\alpha)$, we have
\[ \vol(U_1(\alpha)) = \tau_p(f) p^{(n-1)\alpha}. \]
Now,
\[ U_1(\alpha,\beta)=p^{-\alpha} U_1(0,\beta)=p^{-\alpha}\big(
   \Omega(\beta)-\Omega(\beta+1) \big), \]
therefore
\[ \vol (U_1(\alpha,\beta)) = \frac{p-1}p \tau_p(f) p^{n\alpha-\beta}  . \]

Finally, $U(\alpha)=p^{-\alpha} (\Z_p^n\setminus (p\Z_p^n\cup
\Omega(1)) )$, hence
\[ \vol(U(\alpha))=\big(1-p^{-n} - p^{-1}\tau_p(f)\big)p^{n\alpha} .  \]
\end{proof}

\begin{prop}
Assume that $p\not\in S$. Then,
\[ \hat H_p(\mathbf s;\psi_0) =
\hat H_{\P^n,p}(s_0) + \tau_p(f)
         \frac{p^{s_0-n}-p^{s_1-n}}{(p^{s_0-n}-1)(p^{s_1-n+1}-1)}
\]
where 
\[ \hat H_{\P^n,p}(s_0) = \frac{1-p^{-s_0}}{1-p^{n-s_0}} \]
denotes the Fourier transform (with respect
to the trivial character $\psi_0$)  
of the local height function of $\P^n$
for the tautological line bundle at $s_0$.
\end{prop}
\begin{proof}
By definition, 
\begin{align*}
 \hat H_p(\mathbf s;\psi_0) & = \int_{\Q_p^n} H(\mathbf s;\mathbf x)^{-1}\,
\dx \\
&= \int_{U(0)}  + \sum_{1\leq\beta <\alpha}
 \int_{U_1(\alpha,\beta)} + \sum_{1\leq\alpha} \int_{U_1(\alpha)}
+ \sum_{1\leq\alpha}\int_{U(\alpha)}.
\end{align*}
We compute these sums separately.  The integral over $U(0)$ is equal
to~$1$.
Then 
\begin{align*}
\sum_{1\leq\beta <\alpha}
 \int_{U_1(\alpha,\beta)} &= \frac{p-1}p \tau_p(f)
       \sum_{1\leq\beta<\alpha} 
     p^{-\alpha s_0} p^{-\beta(s_1-s_0)}  p^{\alpha n-\beta} \\
&= \frac{p-1}p \tau_p(f) \sum_{\beta=1}^\infty
     p^{-\beta(s_1-s_0+1)}
 \sum_{\alpha=\beta+1}^\infty p^{-\alpha(s_0-n)} \\
&= \frac{p-1}p \tau_p(f) \sum_{\beta=1}^\infty
     p^{-\beta(s_1-s_0+1)}
     p^{-\beta(s_0-n)} \frac{1}{p^{s_0-n}-1} \\
&= \frac{p-1}p \tau_p(f) \frac{1}{p^{s_0-n}-1}
        \sum_{\beta=1}^\infty p^{-\beta(s_1-n+1)} \\
&= \frac{p-1}p \tau_p(f) \frac{1}{p^{s_0-n}-1} \frac{1}{p^{s_1-n+1}-1}.
\end{align*}
Concerning the integrals over $U_1(\alpha)$, we have
\[
\sum_{1\leq\alpha} \int_{U_1(\alpha)} =
    \tau_p(f) 
        \sum_{\alpha=1}^\infty p^{-\alpha s_1} p^{(n-1)\alpha} 
= \tau_p(f) \frac{1}{p^{s_1-n+1}-1}.
 \]
Finally,
\begin{align*}
\sum_{1\leq\alpha} \int_{U(\alpha)} & =
  (1-p^{-n}-p^{-1}\tau_p(f)) 
\sum_{\alpha=1}^\infty p^{-s_0\alpha} p^{n\alpha}   \\
& = \big( 1-p^{-n}-p^{-1}\tau_p(f) \big) \frac{1}{p^{s_0-n}-1}.
\end{align*}
Adding all these terms gives
\begin{align*}
 \hat H_p(\mathbf s;\psi_0) &= 1 + (1-p^{-1})\tau_p(f)
\frac{1}{p^{s_0-n}-1} \frac{1}{p^{s_1-n+1}-1}  \\
& \qquad{}
+\tau_p(f) \frac{1}{p^{s_1-n+1}-1}
+ \big( 1-p^{-n}-p^{-1}\tau_p(f) \big) \frac{1}{p^{s_0-n}-1} \\
&= 1 + (1-p^{-n}) \frac{1}{p^{s_0-n}-1} \\
& \qquad{}
 + \tau_p(f) \left(
 (1-p^{-1}) \frac{1}{p^{s_0-n}-1} \frac{1}{p^{s_1-n+1}-1}
 + \frac{1}{p^{s_1-n+1}-1} \right. \\
& \qquad\qquad \left. {} - p^{-1} \frac{1}{p^{s_0-n}-1} \right) \\
& = \hat H_{\P^n,p}(s_0) + p^{-1} \tau_p(f)
         \frac{p-1 + p^{s_0-n+1}-p - p^{s_1-n+1}}
         {(p^{s_0-n}-1)(p^{s_1-n+1}-1)} \\
& = \hat H_{\P^n,p}(s_0) + p^{-1} \tau_p(f)
         \frac{p^{s_0-n+1}-p^{s_1-n+1}}{(p^{s_0-n}-1)(p^{s_1-n+1}-1)} \\
& = \hat H_{\P^n,p}(s_0) + \tau_p(f) p^{n-1}
         \frac{p^{-s_1}-p^{-s_0}}{(1-p^{n-s_0})(1-p^{n-1-s_1})}
\end{align*}
\end{proof}

\section{The local Fourier transform at a non-trivial character}

In this subsection  we evaluate the local Fourier transform at 
$p$ for a non-trivial character $\psi_{\mathbf a}$.
Let $S(\mathbf a)$ be the union of $S$ and of the set of primes
$p$ such that $\mathbf a\in p\Z^n$;
We assume that $p\not\in S(\mathbf a)$.

Recall that $Z_f\subset\P^{n-1}_\Z$
denotes the subscheme defined by $f$
and define
$Z_{f,\mathbf a}=Z_f\cap H_{\mathbf a}$,
where $H_{\mathbf a}$ is the hyperplane
of $\P^{n-1}$ defined by $\mathbf a$.
Finally, let $Z_{f,\mathbf a}^t$ (resp.\ $Z_{f,\mathbf a}^{nt}$)
be the locus of points
in $Z_{f,\mathbf a}$
where the intersection $Z_f\cap H_{\mathbf a}$ is \emph{transverse}
(resp.\ is \emph{not transverse}).
By assumption, $Z_f$ and $H_{\mathbf a}$ are smooth over $\Z_p$.

Let $I(\alpha,\beta)$ be the integral of $\psi_{\mathbf a}$
over the set of $\mathbf x\in\Q_p^n$ such
that $\norm{\mathbf x}=p^\alpha$ and
$\abs{f(\mathbf x)}\leq p^{d\alpha -\beta}$.
Then, according to our partition of $\Q_p^n$, we have
\begin{align*}
\hat H_p(\mathbf s;\psi_{\mathbf a} ) & = 1 + \sum_{\alpha=1}^\infty
\sum_{\beta=0}^{\alpha-1} p^{-\alpha s_0} p^{-\beta(s_1-s_0)}
   \int
   _{\substack{\norm{\mathbf x}=p^\alpha \\ 
                \abs{f(\mathbf x)}=p^{d\alpha-\beta}}}
        \psi_{\mathbf a}  \\
&\hskip .5\textwidth
        + \sum_{\alpha=1}^\infty p^{-\alpha s_1} 
        \int_{\substack{\norm{\mathbf x}=p^{\alpha} \\
                \abs{f (\mathbf x)}\leq p^{-\alpha(d-1)}}} 
        \psi_{\mathbf a} \\
&= 1+ \sum_{\alpha=1}^\infty \sum_{\beta=0}^{\alpha-1}
      p^{-\alpha s_0} p^{-\beta(s_1-s_0)}
         \big( I(\alpha,\beta)- I(\alpha,\beta+1)\big) \\
&\hskip .5\textwidth
      + \sum_{\alpha=1}^\infty p^{-\alpha s_1} I(\alpha,\alpha). \\
&= 1+ \sum_{\alpha=1}^\infty p^{-\alpha s_0} I(\alpha,0)
        - (p^{s_1-s_0}-1) \sum_{\alpha=1}^\infty \sum_{\beta=1}^\alpha
             p^{-\alpha s_0} p^{-\beta(s_1-s_0)} I(\alpha,\beta)
\end{align*}

\begin{lem}
If $t\in\Q_p$, the mean value over $\Z_p^*$ of $\psi(t\cdot)$
is equal to
\[ \frac{\int_{\Z_p^*} \psi(tu)\, du}
        {\int_{\Z_p^*} du }
   = \begin{cases}
        1 & \text{if $t\in\Z_p$;} \\
        -1/(p-1) & \text{if $v_p(t)=-1$;} \\
        0 & \text{if $v_p(t)\leq -2$.}
     \end{cases}
\]
\end{lem}
\begin{proof}
Indeed, we have
\begin{align*}
\int_{\Z_p^*} \psi(tu)\, du &=
 \int_{\Z_p} \psi(tu)\, du - \int_{p\Z_p} \psi(tu)\, du \\
&= \int_{\Z_p} \psi(tu)\, du - \frac1p \int_{\Z_p}\psi(ptu)\, du .
\end{align*}
The integral of a non-trivial character over a compact group
is $0$, hence this integral equals
$0$ if $t\not\in p^{-1}\Z_p$, equals $-\frac 1p$ if
$t\in p^{-1}\Z_p\setminus \Z_p$ and equals
$1-\frac1p$ if $t\in\Z_p$. This proves the lemma.
\end{proof}

Using the change of variables $\mathbf x=p^{-\alpha} \mathbf y$,
this implies the following formula:
\begin{multline}
I(\alpha,\beta) = p^{n\alpha} \left(
               \frac p{p-1} \vol\big(\norm{x}=1;\, p^{\beta}|f(\mathbf x);\, 
                p^\alpha|\langle \mathbf a,\mathbf x\rangle \big)\right. \\
\left.  - \frac1{p-1} \vol\big(\norm{x}=1;\, p^{\beta}|f(\mathbf x);\, 
                p^{\alpha-1}|\langle \mathbf a,\mathbf x\rangle \big)
\right).
\end{multline}

\begin{lem}
If $1\leq\beta\leq\alpha$, one has
\[ \vol\big(\norm{\mathbf x}=1;\,
        p^{\beta}| f(\mathbf x);\,
        p^\alpha| \langle\mathbf a,\mathbf x\rangle \big)
 = p^{-\alpha} p^{(2-n)\beta }\big(1-\frac1p\big) \# Z_{f,\mathbf a}(\Z/p^\beta\Z).
\]
\end{lem}

In particular,
 \begin{equation}
 I(\alpha,\beta) = 0 \quad\text{if $1\leq\beta <\alpha$} .
 \end{equation}
Moreover,  if $\alpha\geq 2$,
\begin{align*}
\vol\big( \norm{\mathbf x}=1;\,
        p^{\alpha}| f(\mathbf x);\,
        p^{\alpha-1}| \langle\mathbf a,\mathbf x\rangle\big) 
&= \frac1p \vol\big( \norm{\mathbf x}=1;\,
        p^{\alpha-1}| f(\mathbf x);\,
        p^{\alpha-1}| \langle\mathbf a,\mathbf x\rangle\big) \\
&= \frac1p p^{(1-\alpha)(n-1)} \big(1-\frac1p) \#Z_{f,\mathbf a}(\Z/p^{\alpha-1}) .
\end{align*}
If $\alpha=1$, one has
\[ \vol\big(\norm{\mathbf x}=1;\, 
        p | f(\mathbf x)\big)
     = \big(1-\frac1p\big) p^{1-n} \# Z_f(\Z/p\Z).
\]
We had computed in~\cite{chambert-loir-t99c}, proof of Lemma 3.5, the integral
\[
\int_{\norm{\mathbf x}=p^\alpha}\psi_{\mathbf a}
= \begin{cases} 
        -1       & \text{if $\alpha=1$;} \\
        0        & \text{if $\alpha\geq 2$.}
  \end{cases}
\]
so that
\begin{multline} \label{formula:hathp}
\hat H(\mathbf s;\psi_{\mathbf a}) =
1- p^{-s_0} 
+ \frac{p^{s_1-s_0}-1}{p-1} p^{-s_1} \# Z_f(\F_p) \\
- \frac{p^{s_1-s_0}-1}{p-1} (1-p^{n-s_1-2}) \sum_{\alpha=1}^\infty
        p^{-\alpha(s_1-1)} \# Z_{f,\mathbf a}(\Z/p^\alpha\Z).
\end{multline}

\begin{lem}
For all $\alpha\geq 1$,
\[ \# Z_{f,\mathbf a} (\Z/p^\alpha\Z)
     \leq p^{(n-3)(\alpha-1)} \# Z_{f,\mathbf a}^t (\Z/p\Z)
            + p^{(n-2)(\alpha-1)} \# Z_{f,\mathbf a}^{nt} (\Z/p\Z).
\] 
\end{lem}
\begin{proof}
The inequality is trivially true for $\alpha=1$. We prove it
for any $\alpha$ by induction: to lift a point
in $Z_{f,\mathbf a}(\Z/p^\alpha\Z)$ to a point in
$Z_{f,\mathbf a}(\Z/p^{\alpha+1}\Z)$, one needs to
solve two equations in $\mathbf u\in\F_p^n$:
\[ \langle \nabla f(\mathbf x),\mathbf u\rangle
    \equiv p^{-\alpha} f(\mathbf x) ,
\quad
 \langle {\mathbf a},\mathbf u \rangle
    \equiv  p^{-\alpha} \langle {\mathbf a},\mathbf x\rangle
 \pmod{p} .\]
A point in $Z_{f,\mathbf a}(\Z/p^\alpha\Z)$ which reduces
to a point in $Z_{f,\mathbf a}^t$ modulo $p$ has
$p^{n-3}$ lifts in $Z_{f,\mathbf a}(\Z/p^{\alpha+1}\Z)$.
On the other hand, a point reducing to a point in $Z_{f,\mathbf a}^{nt}$
has $p^{n-2}$ or $0$ lifts according to the two linear equations
being compatible or not. This implies the lemma.
\end{proof}

\begin{prop} \label{prop:37}
If not empty, the set $Z_{f,\mathbf a}^{nt}$ is a closed subscheme of bounded
degree of $Z_{f,\mathbf a}$ and of dimension~$0$.
There exist a constant $C$, independent of $\mathbf a$ and $p$
such that
\[ \# Z_{f,\mathbf a}^t (\Z/p\Z) \leq C p^{n-3},
 \quad \# Z_{f,\mathbf a}^{nt}(\Z/p\Z) \leq C . \]
\end{prop}

As a corollary, one gets:
\begin{cor}
There exist a constant $C$ such that for all $\alpha$ and
$p\not\in S(\mathbf a)$,
\[ \# Z_{f,\mathbf a}(\Z/p^\alpha\Z) 
\leq C p^{(n-3)\alpha} + C p^{(n-2)(\alpha-1)}. \]
\end{cor}

\begin{proof}[Proof of Prop.~\ref{prop:37}]
The set $Z_{f,\mathbf a}$ is defined by the two equations
$f(\mathbf x)=\langle \mathbf a,\mathbf x\rangle=0$.
Fix the coordinates $x_1,\dots,x_n$ so that $\mathbf a$ is the first vector.
Up to a constant, one may write
\[ f(\mathbf x)=x_1^d + g_1(x_2,\dots,x_n)x_1^{d-1} + \cdots + g_{d_1} x_1
+ g_d (x_2,\dots,x_n)
\]
for some homogeneous polynomials $g_i$ of degree $i$.
Then, denoting $\mathbf x= (x_1,\mathbf x')$,
$Z_{f,\mathbf a}$ is defined by the equations
\[ x_1=g_d(\mathbf x')= \partial_2 g_d(\mathbf x') = \dots
=\partial_n g_d(\mathbf x') = 0.
\]
On $Z_{f,\mathbf a}$, $\partial_1 f(0,\mathbf x')=g_{d-1}(\mathbf x')$
and on $Z_{f,\mathbf a}^{nt}\subset Z_{f,\mathbf a}$,
$\partial_i f(0,\mathbf x')=\partial_i g_d(\mathbf x')$. As
$Z_f$ is smooth, $g_{d-1}(\mathbf x')$ doesn't vanish on
$Z_{f,\mathbf a}^{nt}$ which must therefore be either empty
or of dimension $0$.
Its degree cannot exceed $d (d-1)^{n-1}$.
The bound on the number of $\F_p$-rational points are a consequence
of the following (certainly well-known) easy lemma.
\end{proof}
\begin{lem} \label{lem:weil}
Let $k=\F_q$ be a finite field,
$X$ a closed subscheme of $\P^n_k$
of dimension $d$.
Then 
\[ \# X(\F_q) \leq  \P^d(\F_q)  \deg X. \]
\end{lem}
\begin{proof}
We prove this by induction on $d$.
If $d=0$, the result is clear.
Then, one can assume that $X$ is reduced, irreducible
and not contained in any hyperplane.
For any hyperplane $H\subset\P^n$ which is rational over $k$,
$X\cap H$ is a closed subscheme of $H$ of dimension $d-1$ and
of degree $\leq\deg X$.
By induction, we have
\[ \# (X\cap H)(\F_q) \leq \# \P^{d-1} (\F_q) \deg X. \]
Finally, any point of $X(\F_q)$ is contained in exactly $\# \P^{n-1}(\F_q)$
rational hyperplanes in $\P^n$, so that
\[ \# X(\F_q) \#\P^{n-1}(\F_q) \leq \# \P^{d-1}(\F_q) \# \P^n(\F_q) \deg X .\]
As $n\geq d$, this implies
\[ \# X(\F_q)\leq \frac{q^{n+1}-1}{q^n-1} \frac{q^d-1}{q-1} \deg X
\leq \P^{d}(\F_q) \deg X. \]
\end{proof}

\section{The height zeta function}

From now on, we fix some $\eps>0$ and consider only
$\mathbf s$ in the subset $\Omega$ of $\C^2$
defined by the inequalities $\Re(s_0)>n+\eps$ and $\Re(s_1)>n-1+\eps$.

\begin{prop}
There exist a holomorphic function $g$ on $\Omega$ which 
has polynomial growth in vertical strips such that
\[ \hat H(\mathbf s,\psi_0) = g(\mathbf s) \frac{1}{(s_0-n-1)(s_1-n)}.\]
\end{prop}
\begin{proof}
Indeed, we see from 2.3 that for $p\not\in S$,
\[ \hat H_p(\mathbf s,\psi_0) = 
1+ p^{n-s_0} + p^{n-s_1-1} + O(p^{-1-\eps}), \]
the $O$ being uniform in $p$.
Consequently,
\[ \prod_{p\not\in S} \hat H_p(\mathbf s,\psi_0)(1-p^{n-s_0})(1-p^{n-1-s_1})
\]
converges to a holomorphic bounded function on $\Omega$.
As the finite number of remaining factors converge uniformly in $\Omega$,
the existence of $g$ is proven. The growth of $g$ in vertical strips
follows from Rademacher's estimates for the Riemann zeta function.
\end{proof}

\begin{lem}
There exist a constant $C>0$ such that for
all $\mathbf a\in\Z^n\setminus\{0\}$,
all $p\not\in S(\mathbf a)$ and
all $(s_0,s_1)\in\Omega$,
one has
\[ \abs{ \hat H_p(\mathbf s,\psi_{\mathbf a}) -1}
\leq C p^{-1-\eps}. \]
\end{lem}
\begin{proof}
Recall the 
formula~\ref{formula:hathp}:
\begin{multline*}
\hat H_p(\mathbf s,\psi_{\mathbf a}) - 1
=  - p^{-s_0} + \frac{p^{-s_0}-p^{-s_1}}{p-1} p^{n-2}(1-\frac1p)^{-1}\tau_p(f)
\\
- \frac{p^{s_1-s_0}-1}{p-1} (1-p^{n-s_1-2}) \sum_{\alpha=1}^\infty
p^{-\alpha(s_1-1)}\# Z_{f,\mathbf a}(\Z/p^\alpha\Z)
\end{multline*}
the right hand side of which we have to estimate all terms.
The first one is $p^{-s_0}=O(p^{-1-\eps})$.
Then, as $\tau_p(f)$ is bounded, the second one is
\[ O(p^{n-3-\Re(s_0)}) + O(p^{n-3-\Re(s_1)})
      = O(p^{-2}). \]
For the last term $T_3$, we use 
Lemma 3.8 so that, denoting $\sigma_1=\Re(s_1)$,
\begin{multline*}
\sum_{\alpha=1}^\infty
p^{-\alpha(s_1-1)}\# Z_{f,\mathbf a}(\Z/p^\alpha\Z) \\
 \leq C \sum_{\alpha=1}^\infty 
      p^{-\alpha(\sigma_1-1)} p^{(n-3)\alpha}
 + C \sum_{\alpha=1}^\infty 
      p^{-\alpha(\sigma_1-1)} p^{(n-2)(\alpha-1)} \\
\leq C \frac{1}{p^{\sigma_1-n+2}-1} + C p^{2-n} \frac{1}{p^{\sigma_1-n+1}-1}.
\end{multline*}
Moreover, 
\[ \abs{1-p^{n-s_1-2}}  \leq 2  \]
so that
\begin{align*} 
\abs{T_3}&
\ll \frac{1}{p-1}\, \frac{p^{\sigma_1-\sigma_0}+ 1}{p^{\sigma_1-n+2}-1} 
         +2C \frac{p^{2-n}}{p-1} 
\frac{p^{\sigma_1-\sigma_0}+1}{p^{\sigma_1-n+1}-1} \\
&\ll
\frac{1}{p} \big( p^{n-2-\sigma_0} + p^{n-2-\sigma_1}\big)
 + p^{1-n} \big( p^{n-1-\sigma_0}+p^{n-1-\sigma_1}\big) \\
&\ll p^{-2}.
\end{align*}
The lemma is proved.
\end{proof}
 
\begin{prop}
For each $\mathbf a\in\Z^n\setminus\{0\}$,
$\hat H(\mathbf s,\psi_{\mathbf a})$ is
a holomorphic function on $\Omega$.
Moreover, there exist constants $C>0$  and $\nu$
(which are independent of $\mathbf s$ and $\mathbf a$) such that
\[ \abs{\hat H(\mathbf s,\psi_{\mathbf a})}
        \leq C (1+\norm{\Im(s)})^\nu (1+\norm{\mathbf a})^{-n-1}.\]
\end{prop}
\begin{proof}
Write
\[ \hat H(\mathbf s,\psi_{\mathbf a}) 
= \prod_{p\not\in S(\mathbf a)} \hat H_p
   \times 
    \prod_{p\in S(\mathbf a)} \hat H_p
     \times
  \hat H_\infty.
\]
The convergence of the first infinite product to a bounded
holomorphic function follows from the preceding
lemma. 
As in Lemma~3.7 of \cite{chambert-loir-t99c}, there exists
a constant $\kappa>0$ such that
\[ \abs{\prod_{p\in S(\mathbf a)} \hat H_p (\mathbf s,\psi_{\mathbf a})}
\ll (1+\norm{\mathbf a})^{\kappa}.
\]
Using the rapidly decreasing
behaviour of $\hat H_\infty$ as a function of $\mathbf a$
\[ \abs{\hat H_\infty(\mathbf s,\psi_{\mathbf a})}
         \ll (1+\norm{\mathbf a})^{-n-\kappa-1} \]
established in Prop.~2.13 of \emph{loc. cit.}, the proposition
is proved.
\end{proof}

\begin{thm} \label{thm:peyre}
The height zeta function converges in the domain $\Re(s_0)>n+1$,
$\Re(s_1)>n$. Moreover, there
exists a holomorphic function $g$ in the domain $\Re(s_0)>n$,
$\Re(s_1)>n-1$ such that
\[ Z(\mathbf s) = g(\mathbf s) \frac{1}{(s_0-n-1)(s_1-n)}. \]
The function $g$ has polynomial growth in vertical strips
and $g(n+1,n)\neq 0$.
\end{thm}

Specializing to $\mathbf s=s(n+1,n)$ and using a standard Tauberian
theorem, one obtains the following corollary.
\begin{cor} \label{cor:peyre}
There exist a polynomial $P_X$ of degree $1$ and a real number
$\alpha>0$ such that the number of points of 
$U(\Q)\subset X(\Q)$
of anticanonical height $\leq B$ satisfies
\[ N(U,K_X^{-1},B)=B P(\log B) + O(B^{1-\alpha}). \]
Moreover, if $\tau(K_X)$ denotes the Tamagawa number, the
leading coefficient of $P_X$ is equal to
\[ \frac{\tau(K_X)}{(n+1)n}, \]
as predicted by Peyre's refinement of Manin's conjecture.
\end{cor}

\def\noop#1{\ignorespaces}
\bibliographystyle{smfplain}
\bibliography{acl}
\end{document}